# On the Mathematical Modelling of Measurement

Jonathan Barzilai[1]

## Abstract

The operations of linear algebra, calculus, and statistics are routinely applied to measurement scales but certain mathematical conditions must be satisfied in order for these operations to be applicable. We call attention to the conditions that lead to construction of measurement scales that enable these operations.

## 1 Introduction

The problem of applicability of mathematical operations to scale values has received little attention in the literature. The purpose of this paper is to summarize the conditions for the construction of measurement scales that enable the application of the operations of algebra, calculus, and statistics to scale values. For additional details in the context of decision theory and its applications see Barzilai [1–3].

## 2 The Purpose of Measurement

Some of the universally accepted key elements of the Theory of Measurement have been identified as early as 1887 by Helmholtz [5] and 1920 by Campbell [4], but the very first section of von Neumann and Morgenstern [7, §1] leaves no doubt that difficulties with the problem of mathematical modelling of measurement remained.

To clarify what is meant by "the mathematical modelling of measurement" some terminology is required. By an empirical system $E$ we mean a set of empirical objects together with operations (i.e. functions) and possibly the relation of order (which is not an operation since it is not a single-valued function). A mathematical model $M$ of the empirical system $E$ is a set with operations that reflect the operations in $E$ as well as the order in $E$ when $E$ is ordered. A scale $s$ is a mapping of the objects in $E$ into the objects in $M$ that reflects the structure of $E$ into $M$ (technically, a scale $s$ is a homomorphism from $E$ into $M$) and we say that $s$ preserves the structure of $E$. The purpose of this construction is to enable the application of mathematical operations on the elements of the mathematical system $M$ (which, for the reasons given in §4, typically is the set of real numbers with appropriate operations).

As Campbell eloquently states [4, pp. 267–268], "the object of measurement is to enable the powerful weapon of mathematical analysis to be applied to the subject matter of science." Although these concepts have been universally accepted, and despite the fact that "the powerful weapon of mathematical analysis" cannot be applied if mathematical operations are not applicable in the system $M$, fundamental questions related to applicability of mathematical operations have received little attention in the literature. Consider the applicability of the operations of addition and multiplication on scale values for a fixed scale, that is, operations that express facts such as "the weight of a given object equals the sum of the weights of two other objects" ($s(a) = s(b) + s(c)$) and "the weight of a given object is two and a half times the weight of another one" ($s(a) = 2.5 \times s(b)$). Models that do not enable the operations of addition and multiplication cannot serve as the foundation of any scientific discipline since much of mathematical analysis is not applicable when linear algebra and, *a fortiori*, calculus, are not applicable.

---

1  Dept. of Industrial Engineering, Dalhousie University, Halifax, Nova Scotia, B3J 2X4 Canada.
   © J. Barzilai 2006. Email: Barzilai@dal.ca

It is important to emphasize the distinction between the application of the operations of addition and multiplication on scale values for a fixed scale (for example $s(a) = s(b) + s(c)$) as opposed to what appear to be the same operations when they are applied to an entire scale whereby an equivalent scale is produced by what amounts to a change of the zero point or unit (for example $t = p + q \times s$ where $s$ and $t$ are two scales and $p, q$ are numbers). In the first case, the operations of addition and multiplication are applied to elements of the system $M$ and the result is another element of $M$. In the latter case, addition or multiplication by a number are applied to an element of the set $S = \{s, t, …\}$ of all possible scales and the result is another element of $S$ rather than $M$. Because operations are functions, and functions with different domains or ranges are different, these are different operations. In the case of "interval" scales where the uniqueness of the set of all possible scales is characterized by $t = p + q \times s$, it cannot be concluded that the operations of addition and multiplication are applicable on scale values for a fixed scale such as $s(a) = s(b) + s(c)$. It might be claimed that the characterization of scale uniqueness by $t = p + q \times s$ implies the applicability of addition and multiplication on scale values for fixed scales, but this claim requires proof. (There is no such claim, nor such proof, in the literature because a simple argument[2] shows that this claim is false.)

## 3    The Principle of Reflection

Recall that the purpose of measurement is to enable the application of mathematical operations in the system $M$. If the operations of addition and multiplication are to be applicable in $M$ in their usual algebraic form, including subtraction and division, $M$ must be equipped with these operations. The technical name for such a system is *field*. In simple terms, $M$ behaves like a two-button calculator: one for addition, the other for multiplication. The significance of this observation is that empirical systems that are based on a single operation (one-button calculators) cannot be reflected, i.e. modelled, by mathematical systems with two operations (two-button calculators). In other words, if two operations are reflected in *M*, there must be two operations in *E* to be reflected. For a detailed discussion see *The Principle of Reflection* in Barzilai [1, §6.8]. The unavoidable conclusion is that models of decision (and measurement) theory that are based on any single operation do not enable the application of the operations of addition and multiplication. (In technical terms, a field cannot be the homomorphic image of a group).

---

2 Consider the automorphisms of the group of integers under addition. The group is a model of itself ($E = M$), and scale transformations are multiplicative: $t = (\pm 1) \times s$. However, by definition, the operation of multiplication which is defined on the set of scales is not defined on the group $M$.

It may be argued that it is sufficient to model the operation of addition since multiplication is repeated addition, but this is only true for the natural numbers. In general, and in particular for the real numbers, multiplication is not defined as repeated addition but through field axioms. In simple terms, either *M* is a one-button calculator or it is a two-button calculator, but it cannot be both.

## 4 Homogeneity Considerations

When the empirical system *E* is ordered, its model *M* must be ordered as well so that the underlying field is ordered. Some fields are ordered, (e.g. the field of rational numbers) while others (for example the field of complex numbers) are not. In addition, in order for "the powerful weapon" of calculus to be applicable in *M*, the limit operation must be enabled – addition, multiplication, and order are not sufficient to enable the application of the limit operation in the mathematical system *M* as can be seen from the example of the field of rational numbers. In technical terms, the limit operation requires that the underlying field be *complete* and since the only ordered, complete field is the field of real numbers (see e.g. McShane and Botts [6, pp. 22–24]), we conclude that the application of calculus in an ordered field can only be carried out in the field of real numbers.

Finally, when homogeneity considerations (see below – these considerations concern the existence of an absolute zero and multiplicative unit in *E* and its model *M)* are taken into account, it turns out (see Barzilai [1, §6.10 and 2, pp. 176–177]) that there are three models that enable the application of calculus: The measured objects must correspond to (i) scalars in the field of real numbers; or to (ii) vectors in a one-dimensional vector space over this field; or to (iii) points in a one-dimensional affine space over the real numbers. (Technically, the zero vector in a vector space is an absolute zero because it is a fixed point of the automorphisms of the space.)

For psychological variables where the existence of an absolute zero is not established, the only possibility for addition, multiplication, order, and limits to be applicable is the model where the measured objects correspond to points in a one-dimensional affine space over the ordered real numbers. In such a space the ratio of two points is undefined while their difference is a vector and the ratio of two vectors is a scalar.

Ratios of variables for which the existence of an absolute zero has not been established are undefined. For example, ratios $T_1/T_2$ of temperature have been undefined until it was established that temperature has an absolute zero (see e.g. von Neumann and Morgenstern [7, §3.4.6]). In the case of *time*, the ratio $t_1/t_2$, where $t_1$ and $t_2$ are two *points* in time, is undefined while the ratio $(\Delta t)_1/(\Delta t)_2$ of two time *differences*, i.e. time periods or time intervals, is well-defined. It follows that the ratio $v_1/v_2$ is undefined for any psychological variable since the existence of an absolute zero has not been established for psychological variables.

## 5  Conclusions

In the case of physical variables, the set of scales is uniquely determined by the set of objects and the property under measurement. In other words, scale construction requires specifying only the set of objects and the property under measurement. In the social sciences, the system under measurement includes a person or persons so that the property under measurement is associated with a human being and, in this sense, is personal, psychological, or subjective. Except that in the case of subjective properties the specification of the property under measurement includes the specification of the "owner" of the property (for example, we must specify whose preference is being measured), the mathematical modelling of measurement of subjective properties does not differ from that of physical ones.

In both cases, i.e. for physical as well as psychological variables (including preference), a necessary condition for the application of linear algebra and calculus to the values of measurement scales is that addition and multiplication be defined in the empirical system. For variables which do not have an absolute zero or unit such as *preference, time, potential energy, position,* the correct model is a one-dimensional affine space (see Barzilai [1–3] for details). In the empirical systems and their mathematical models for affine variables, *ratios* are undefined but *differences* and *ratios of differences* of affine variables are well-defined.